\documentclass{article}
\usepackage{amsfonts,amsmath,amssymb,amsthm,cite,enumerate,graphicx,mathrsfs}
\usepackage{authblk}
\usepackage{chngcntr}
\counterwithin{figure}{section}
\usepackage{float}
\usepackage{geometry}
\date{}
\makeatletter

\newcommand{\Rmnum}[1]{\expandafter\@slowromancap\romannumeral #1@}
\makeatother
\usepackage{lmodern}
\usepackage[utf8]{inputenc}
\usepackage{hyperref}

\numberwithin{equation}{section}
\newtheorem{corollary}{Corollary}[section]

\newtheorem{lemma}{Lemma}[section]
\newtheorem{proposition}{Proposition}[section]
\newtheorem{remark}{Remark}[section]

\newtheorem{theorem}{Theorem}[section]

\newcommand\theref[1]{Theorem~\ref{#1}}
\newcommand\lemref[1]{Lemma~\ref{#1}}

\title{Perturbation Analysis of Gravitational Collapse of Gaseous Stars in the Pressureless Euler-Poisson System}
\author{Xingchen Wang\thanks{{e}-mail: wxc3619@163.com} \quad Xiaomin Zhang\thanks{Corresponding author {e}-mail: zxm15924687@163.com}
 \\ \small\textit{ Department of mathematics, Shandong Normal University, Jinan 250014 China}}
\begin{document}
\begin{sloppypar}
\date{}
\maketitle
\begin{center}
\begin{minipage}{130mm}{\small
\textbf{Abstract}:
The compressible Euler-Poisson system describes the full evolutionary processes of inviscid compressible gases, including compression, expansion, fluid motion, equilibrium states and gravitational collapse, under the combined effects of self-gravitation, pressure and inertial forces. It constitutes the fundamental macroscopic mathematical model for describing stellar structures, nebular evolution and gravitational collapse of gaseous celestial bodies. The system comprises mass conservation, momentum conservation, energy conservation laws coupled with the Poisson equation. Following the research in \cite{huang2022}, this paper investigates the quantitative relations among the relative perturbation amplitude of initial velocity, critical time and collapse time for gaseous stars undergoing gravitational collapse within the framework of the pressureless Euler-Poisson system. The closed-form formulas are verified against direct numerical integration, the accuracy of the asymptotic power law is quantified, and the relation between the two time scales is examined. 
}
\end{minipage}
\end{center}

\section{Introduction}
\indent\indent Gravitational collapse of stars is one of the most spectacular and pivotal phenomena in astrophysics, representing the terminal evolutionary stage of massive stars. Once nuclear fusion fuel in the stellar core is depleted, radiative pressure can no longer balance the star's self-gravity, which initiates rapid collapse of the stellar core. Studies on the dynamics of gravitational collapse are fundamentally important for revealing the formation of cosmic structures and characterizing material behaviors under extreme physical conditions.

This paper is concerned with the three-dimensional compressible Euler-Poisson system
\begin{equation}\label{1.1}
\left\{
\begin{array}{ll}
\partial _{t}\rho +\nabla (\rho \mathbf{u})=0,\\
(\rho u)_{t}+\nabla \cdot (\rho u\otimes u)+\nabla P+\rho \nabla \Phi =0,\\
\Delta \Phi =4\pi \rho , \lim\limits_{|\mathbf{x}|\to \infty }\Phi (t,\mathbf{x})=0,
\end{array}
\right.
\end{equation}
where, $\rho$ denotes the mass density, $\mathbf{u}=(u_{1}, u_{2}, u_{3})$ stands for the fluid velocity vector, $\Phi$ represents the gravitational potential, $\Delta=\nabla^{2}$ denotes the Laplace operator, $P$ is the fluid pressure obeying a polytropic equation of state
\begin{equation}\label{1.2}
P=P(\rho)=\rho^{\gamma}, 
\end{equation}
where $\gamma>1$ denotes the adiabatic index.

This system combines the Euler equations from fluid mechanics with the Poisson equation describing the gravitational field, effectively modeling the density evolution, velocity distribution, and changes in the gravitational field of self-gravitating fluids, and has significant applications in astrophysics and theoretical physics.
Early research focused on modeling stellar equilibrium states, and investigations into classical solutions of this system originated from equilibrium problems of gaseous stars. Specifically, non-rotating spherically symmetric static solutions satisfy the Lane-Emden equation, which forms an essential theoretical foundation for understanding the internal structure of stars. On this basis, Chandrasekhar \cite{chandrasekhar1967} systematically elaborated the fundamental theory of stellar structure and incorporated Lane-Emden theory into a rigorous astrophysical framework.

In the development of the theoretical framework for rotating stars, Auchmuty and Beals \cite{auchmuty1971} adopted variational methods to study a class of nonlinear free boundary problems and established an existence theory for equilibrium solutions of rotating self-gravitating fluids. Caffarelli and Friedman \cite{caffarelli1980} further conducted geometric and regularity analyses of free boundaries for axisymmetric rotating self-gravitating fluids and proved under suitable technical assumptions that equilibrium configurations contain at most finitely many fluid rings whose free-boundary cross-sections are analytic on each ring's interior yet may develop singular cusps at ring endpoints. From the perspective of astrophysical modeling and macroscopic dynamics, Binney and Tremaine \cite{binney2008} systematically summarized various steady-state dynamical models of rotating galaxies and constructed a theoretical framework for self-gravitating rotating celestial bodies. For refined analysis of uniformly rotating stellar models, Chanillo and Li \cite{chanillo1994} further explored the geometric scales and morphological features of uniformly rotating stars and derived a priori estimates for stellar diameters.

Stability theory represents an important research branch concerning the Euler-Poisson system. Building on \cite{deng2002}, which addressed the existence and uniqueness of classical solutions for isentropic gaseous star Euler-Poisson equations, Xie and Li \cite{xie2012} extended the work to steady-state solutions of the non-isentropic case. Jang \cite{jang2014} studied the nonlinear instability of Lane-Emden stars, uncovered the mechanisms responsible for stellar structural instability under specific parameter regimes.

The physical vacuum problem presents a major research challenge within this field. Coutand and Shkoller \cite{coutand2012} achieved remarkable progress in the physical vacuum problem for the three-dimensional compressible Euler equations with moving boundaries and established well-posedness within suitable smooth function spaces. Sideris \cite{sideris2017} studied the global existence and asymptotic behavior of affine motions of three-dimensional ideal fluids surrounded by vacuum, and obtained results describing the long-time evolution of fluids subject to vacuum boundary conditions.

Beyond the well-posedness and long-time evolution at physical vacuum boundaries, which are widely recognized as formidable open problems in this field, the polytropic exponent $\gamma$ serves as a core parameter of the Euler-Poisson system governing the dynamical evolution of gaseous stars. Under the adiabatic approximation, the polytropic exponent coincides with the adiabatic index of the gas and its value fundamentally alters three key properties of solutions for the system: the blowup mechanism, the existence of self-similar structures, and density concentration behavior. Goldreich and Weber \cite{goldreich1980} demonstrated that homologous self-similar collapse solutions exist for the adiabatic index $\gamma=\frac{4}{3}$ in the equation \eqref{1.2}, which acts as the critical threshold for the stability of the system. Comprehensive investigations on this family of solutions have been carried out in \cite{deng2003, makino1992, fu1998}. Hadžić and Jang \cite{hadzic2018} classified the adiabatic index $\gamma$ into three regimes: the subcritical case ($\gamma>\frac{4}{3}$), the critical case ($\gamma=\frac{4}{3}$), and the supercritical case ($1<\gamma<\frac{4}{3}$). Guo, Hadžić, and Jang \cite{guo2020} made pivotal contributions to the dynamical research concerning the supercritical regime $1<\gamma<\frac{4}{3}$. They systematically investigated the continued gravitational collapse of gaseous stars governed by the pressureless Euler-Poisson system.
Building upon the work \cite{guo2020}, Huang and Yao \cite{huang2022} conducted more general and complicated mathematical investigations on the continued gravitational collapse of gaseous stars within the framework of the pressureless Euler-Poisson system. The authors identify the physical mechanism for sustained gravitational collapse of gaseous stars with \(1<\gamma<\frac{4}{3}\). They explicitly construct an infinite-dimensional family of collapse solutions to the pressureless Euler-Poisson equations and classify all spherically symmetric solutions into escape and collapse cases.

Although considerable progress has been made in the solution theory, stability analysis, and critical collapse mechanisms of the Euler-Poisson system, several essential research gaps remain unaddressed. Current literature fails to explicitly characterize how initial flow configurations quantitatively affect the critical blow-up time during gravitational collapse of gaseous stars. Furthermore, systematic numerical frameworks and quantitative prediction tools dedicated to such collapse dynamics are still lacking. These limitations constitute a key bottleneck that restricts further refinement of the theoretical framework for stellar gravitational collapse.

Accordingly, building on the theoretical framework proposed in \cite{huang2022}, this paper investigates the pressureless Euler-Poisson system in the supercritical regime $1<\gamma<\frac{4}{3}$ via a combination of theoretical derivation and numerical validation. We mainly focus on gravitational collapse problems with initial velocity fields located within the critical parameter range. By defining the relative perturbation amplitude of the initial velocity relative to the critical velocity threshold --- hereafter termed the perturbation amplitude parameter --- this work quantifies the systemic deviation from equilibrium and establishes a valid stability criterion. Furthermore, explicit quantitative correlations among perturbation amplitude, critical blow-up time, and gravitational collapse time are derived. Combined with a numerical verification of the closed-form formulas, the accuracy of the asymptotic law is assessed. This study further improves the theoretical system of gaseous star collapse governed by the pressureless Euler-Poisson model and provides innovative analytical insights and quantitative approaches for future relevant investigations.

The paper is organized as follows: Section 2 presents preliminary analyses of gaseous star collapse based on the pressureless Euler-Poisson system. Section 3 establishes explicit quantitative relations among the perturbation amplitude parameter, critical time and collapse time. Section 4 verifies the closed-form formulas numerically, quantifies the accuracy of the asymptotic law and compares the two time scales.

\section{Preliminaries}
\subsection{The supercritical regime of gaseous stars}

\indent\indent  We define $\Omega(t)$ as the support of the density function $\rho(t, \cdot)$, with $\partial \Omega(t)$ representing its boundary, and $n(t)$ denoting the unit outward normal vector on $\partial \Omega(t)$. Then, the physical vacuum condition is formulated as
\[
-\infty<\nabla\left(\frac{d p}{d \rho}(\rho)\right) \cdot n\big|_{\partial \Omega(t)}<0,
\]
along with the boundary decay behavior of density
\[
\lim _{\Omega(t) \to \partial \Omega(t)} \rho=0.
\]
This set of boundary conditions describes the smooth transition of density from the stellar surface to the surrounding vacuum. To derive the pressureless Euler-Poisson system in which the pressure contribution is negligible, as \cite{huang2022}, we introduce a radial scaling factor $\mathcal{A}(s, r)$ subject to the initial data
\[
\begin{cases}
\mathcal{A}(0, r) & =\mathcal{A}_{0}(r)>0, \\
\mathcal{A}_{s}(0, r) & =\mathcal{A}_{1}(r),
\end{cases}
\]
and construct a special coordinate mapping
\[
\tau (s,y)=\mathcal{A} (s,r)y, 
\]
with $r=|y|,\;r\in [0,1]$. Under the Lagrangian coordinate transformation, the Jacobian determinant of the mapping is computed as
\[
J[\mathcal{A}]=det(D \tau)=\mathcal{A}^{2}\left(\mathcal{A}+r \partial_{r} \mathcal{A}\right).
\]
Fix an arbitrary small parameter $\bar{\varepsilon}>0$ and set $\varepsilon=\bar{\varepsilon}^{4-3 \gamma}$. After implementing dimensionless scaling transformations, we obtain the pressureless Euler-Poisson system formulated in Lagrangian coordinates
\begin{equation}\label{2.1}
\left\{
\begin{array}{ll}
\mathcal{A} _{ss}+\frac {G(r)}{\mathcal{A} ^{2}}=0,\\
\Delta \tilde {\Phi }=4\pi \tilde {\rho },\\
\lim\limits_{|x|\to \infty }\tilde {\Phi } (t,x)=0,
\end{array}
\right. 
\end{equation}
where the gravitational potential is rescaled as $\tilde{\Phi}=\tilde{\Phi}(s,y)=\bar{\varepsilon}\Phi(t,x)$ with $s=\bar{\varepsilon}^{-\frac{3}{2}}t$, $y=\bar{\varepsilon}^{-1}x$, the mass density satisfies $\tilde{\rho}=\bar{\varepsilon}^{3} \rho(s, y)$, the pressure functional is given by
\[
p[\mathcal{A}]:=\frac{\mathcal{A}^{2}}{\omega^{\alpha} r^{2}}\left(r \partial_{r}\right)\left(\omega^{1+\alpha} J[\mathcal{A}]^{-\gamma}\right),
\]
and
\begin{equation}\label{2.2} 
G(r) \triangleq \frac{1}{r^{3}} \int_{0}^{r} 4 \pi \omega^{\alpha} \delta^{2} d\delta
\end{equation}
denotes the mean density enclosed within the radial ball of radius $r$, where $\omega(r)$ stands for the fluid enthalpy and the exponent satisfies $\alpha=\frac{1}{\gamma-1}$.

To guarantee the physical validity and well-posed solvability of the pressureless Euler-Poisson system, appropriate initial conditions must be imposed. As discussed above, we prescribe the initial data
\begin{equation}\label{2.3}
\left\{
\begin{array}{l}
\mathcal{A}(0, r)=\mathcal{A} _{0}(r)>0, \\
\mathcal{A}_{s}(0, r)=\mathcal{A} _{1}(r), \\
\tilde {\rho }(0, r)=\tilde {\rho }_{0}(r).
\end{array}
\right.
\end{equation}

Owing to energy conservation, substituting the initial energy identity
\begin{equation}\label{2.4}
E(0)=\frac{1}{2} \mathcal{A}_{1}^{2}-\frac{G(r)}{\mathcal{A}_{0}}, 
\end{equation}
yields the critical velocity
\begin{equation}\label{2.5}
v_{0}(r)=\sqrt{2 G(r)}. 
\end{equation}
During the gravitational collapse of gaseous stars, equation \eqref{2.5} acts as the core criterion to judge whether stellar matter at the corresponding radius can sustain continuous collapse.

Moreover, it follows from the law of energy conservation that
\begin{equation}\label{2.6}
\mathcal{A}_{s}^{2}=\mathcal{A}_{1}^{2}+2 G(r)\left(\frac{1}{\mathcal{A}}-\frac{1}{\mathcal{A}_{0}}\right).
\end{equation}
For simplicity, we set $\mathcal{A}_{0}(r)=1$ throughout this paper. Then, the particle velocity evolution equation can be readily derived as
\begin{equation}\label{2.7}
\mathcal{A}_{s}(s, r)=\pm \sqrt{\mathcal{A}_{1}^{2}(r)-2 G(r)+\frac{2 G(r)}{\mathcal{A}(s, r)}}.
\end{equation}
In addition, we denote the initial energy parameter as $n_{0}$, namely
\begin{equation}\label{2.8}
n_{0}=\mathcal{A}_{1}^{2}(r)-2 G(r). 
\end{equation}

When the initial energy parameter satisfies $n_{0}<0$ and the initial velocity fulfills $0<\mathcal{A}_{1}(r)<v_{0}(r)$, the system first expands until reaching the critical time $s_{1}(r)$. At this instant, $\mathcal{A}_{s}(s_{1}, r)=0$, and the system attains its maximum expansion state. Subsequently, the system begins to contract and eventually undergoes complete gravitational collapse at time $s_{0}(r)$, where $\mathcal{A}(s_{0}, r)=0$. The entire dynamical evolution is dominated by the perturbation amplitude parameter $\xi(r)=1-\frac{\mathcal{A}_{1}(r)}{v_{0}(r)}$, whose physical significance will be elaborated in detail in the subsequent Section.

\subsection{Sustainability of gravitational collapse}
\indent\indent For spherically symmetric solutions of the pressureless Euler-Poisson system, the monotonicity of the collapse time $s_{0}(r)$ is a fundamental property for the analysis of gravitational collapse dynamics. By virtue of proof by contradiction, \cite{huang2022} rigorously demonstrates that the collapse time $s_{0}(r)$ is non-decreasing with respect to the radial variable $r$, which indicates the persistence of collapse solutions for the pressureless Euler-Poisson system.

Inspired by the work in \cite{huang2022}, we point out that non-collision of trajectories serves as a necessary condition for the progression of gravitational collapse in the study of spherically symmetric solutions to the pressureless Euler-Poisson system. During the gravitational collapse of gaseous stars, the trajectories of particles labelled by distinct Lagrangian coordinates must be strictly monotonic in physical space to guarantee the physical validity of the system.

\begin{lemma}\label{lem:2.1}
In Lagrangian coordinates, the averaged mass function with respect to $r$ satisfies
\begin{equation}\label{2.9}
G'(r)<0.
\end{equation}
\end{lemma}

\begin{proof}
 Let $\theta=\frac{\delta}{r}$, it yields from \eqref{2.2}
\[
\begin{aligned}
G(r) & =\frac{1}{r^{3}} \int_{0}^{1} 4 \pi \omega^{\alpha}(r \theta)(r \theta)^{2} \cdot r d \theta \\
& =\int_{0}^{1} 4 \pi \omega^{\alpha}(r \theta) \theta^{2} d \theta .
\end{aligned}
\]
Differentiating $G(r)$ with respect to $r$, we obtain
\[
\begin{aligned}
G'(r) & =\frac{d}{d r} \int_{0}^{1} 4 \pi \omega^{\alpha}(r \theta) \theta^{2} d \theta \\
& =\int_{0}^{1} 4 \pi \frac{d \omega^{\alpha}}{d r} \theta^{3} d \theta.
\end{aligned}
\]
We assume that
\[
\frac{d \omega^{\alpha}}{d r}<0, \quad r \in(0,1],
\]
which leads to
\[
G'(r)<0.
\]
Thus \lemref{lem:2.1} is proved.
\end{proof}

\begin{theorem}\label{thm:2.1}
The sufficient condition ensuring non-intersecting particle trajectories reads
\begin{equation}\label{2.10}
\frac{d \mathcal{A}_{1}(r)}{d r}>0.
\end{equation}
\end{theorem}

\begin{proof}We proceed by contradiction. Consider two distinct initial radii satisfying $r_{1}<r_{2}$, whose corresponding trajectories are $\mathcal{A}(s, r_{1})$ (inner layer) and $\mathcal{A}(s, r_{2})$ (outer layer). If a collision occurs, there exists some time $s_{c}$ such that
\[
\mathcal{A} (s_{c},r_{1})=\mathcal{A} (s_{c},r_{2}).
\]
Since $\mathcal{A}(0, r)=1$, we perform a Taylor expansion for small $s$ to obtain $\mathcal{A}(s, r)=1+\mathcal{A}_{1}(r) s+o(s^{2})$. For $r_{1}<r_{2}$, it follows that $\mathcal{A}(s, r_{2})-\mathcal{A}(s, r_{1})=[\mathcal{A}_{1}(r_{2})-\mathcal{A}_{1}(r_{1})] s$. From
$\frac{d\mathcal{A}_{1}(r)}{d r}= \lim\limits_{r_2-r_1 \to 0} \frac{\mathcal{A}_{1}(r_{2})-\mathcal{A}_{1}(r_{1})}{r_{2}-r_{1}}>0$, we have
\begin{equation}\label{2.11}
\mathcal{A}_{1}(r_{2})>\mathcal{A}_{1}(r_{1}).
\end{equation}
Then for all $s<s_{c}$, we have $\mathcal{A}(s, r_{1})<\mathcal{A}(s, r_{2})$. Furthermore at the collision instant $s_{c}$, the inequality
\begin{equation}\label{2.12}
\frac{\mathcal{A}(s_{c}, r_{1})-\mathcal{A}\left(s, r_{1}\right)}{s_{c}-s}>\frac{\mathcal{A}\left(s_{c}, r_{2}\right)-\mathcal{A}\left(s, r_{2}\right)}{s_{c}-s}
\end{equation}
holds, which implies
\begin{equation}\label{2.13}
\mathcal{A} _{s}(s_{c},r_{1})>\mathcal{A} _{s}(s_{c},r_{2}).
\end{equation}
Combining \eqref{2.11} and \eqref{2.13} yields
\begin{equation}\label{2.14}
\mathcal{A} _{s}(s_{c},r_{1})-\mathcal{A} _{1}(r_{1})>\mathcal{A} _{s}(s_{c},r_{2})-\mathcal{A} _{1}(r_{2}).
\end{equation}
By \eqref{2.9}, it yields $G(r_{1})>G(r_{2})$ when $r_{1}<r_{2}$. Then for any $s<s_{c}$, by $\eqref{2.1}_{1}$, we get
\begin{equation}\label{2.15}
\mathcal{A}_{s s}(s, r_{1})=-\frac{G\left(r_{1}\right)}{\mathcal{A}^{2}\left(s, r_{1}\right)}<-\frac{G\left(r_{2}\right)}{\mathcal{A}^{2}\left(s, r_{2}\right)}=\mathcal{A}_{s s}\left(s, r_{2}\right).
\end{equation}
Furthermore, it yields
\begin{equation}\label{2.16}
\left[\mathcal{A} _{s}(s_{c},r_{1})-\mathcal{A} _{1}(r_{1})\right]-\left[\mathcal{A} _{s}(s_{c},r_{2})-\mathcal{A} _{1}(r_{2})\right]=\int _{0}^{s_{c}}\mathcal{A} _{ss}(s,r_{1})ds-\int _{0}^{s_{c}}\mathcal{A} _{ss}(s,r_{2})ds<0,
\end{equation}
namely
\begin{equation}\label{2.17}
\mathcal{A} _{s}(s_{c},r_{1})-\mathcal{A} _{1}(r_{1})<\mathcal{A} _{s}(s_{c},r_{2})-\mathcal{A} _{1}(r_{2}),
\end{equation}
which contradicts \eqref{2.14}. Therefore, the initial assumption is false, and no particle collision can take place. This completes the proof of \theref{thm:2.1}.
\end{proof}

This condition ensures that particles at outer Lagrangian radii always lie outside those at inner radii, thereby preserving spherical symmetry. When the initial velocity $\mathcal{A}_{1}(r)$ satisfies \eqref{2.10}, particle collisions are excluded, and persistent gravitational collapse proceeds instead. This guarantees that the gaseous star maintains a well-behaved geometric configuration throughout collapse, which is of great significance for understanding the intrinsic nature of gravitational collapse.

\section{Quantitative relations among perturbation amplitude parameter, critical time and collapse time}
\subsection{Definition and characterization of the perturbation amplitude parameter}
\indent\indent When investigating the gravitational collapse of gaseous stars, the degree to which the system deviates from equilibrium constitutes one of the key factors governing the dynamical collapse process. To quantify this deviation, we define the perturbation amplitude parameter $\xi(r)$, which acts as the dominant perturbation source driving the transition from initial expansion to collapse. Physically, it represents a normalized velocity deficit, given by
\begin{equation}\label{3.1}
\xi(r)=1-\frac{\mathcal{A}_{1}(r)}{v_{0}(r)}.
\end{equation}
The magnitude of the perturbation amplitude parameter $\xi(r)$ is subject to constraints. Based on the stability analysis of the spherically symmetric Euler-Poisson system, the perturbation amplitude parameter can only take values within a finite interval, namely
\begin{equation}\label{3.2}
0<\xi(r)<1.
\end{equation}
When $\xi(r)=0$, the system is in a critical escape state with no collapse occurring. As $\xi(r)$ increases, the deviation of the system from equilibrium intensifies, and the gravity-induced collapse effect becomes stronger. When $\xi(r)=1$, the system loses stability and pure gravitational collapse takes place.

The perturbation amplitude parameter serves to quantify the decisive influence of initial conditions on the evolutionary path of gaseous stars. Small perturbations may trigger stellar oscillations, while sufficiently large perturbations inevitably lead to gravitational collapse of the star. This parameter acts as a critical criterion for judging the long-term dynamical fate of the star.

\begin{figure}[H]
\centering
\includegraphics[width=0.82\textwidth]{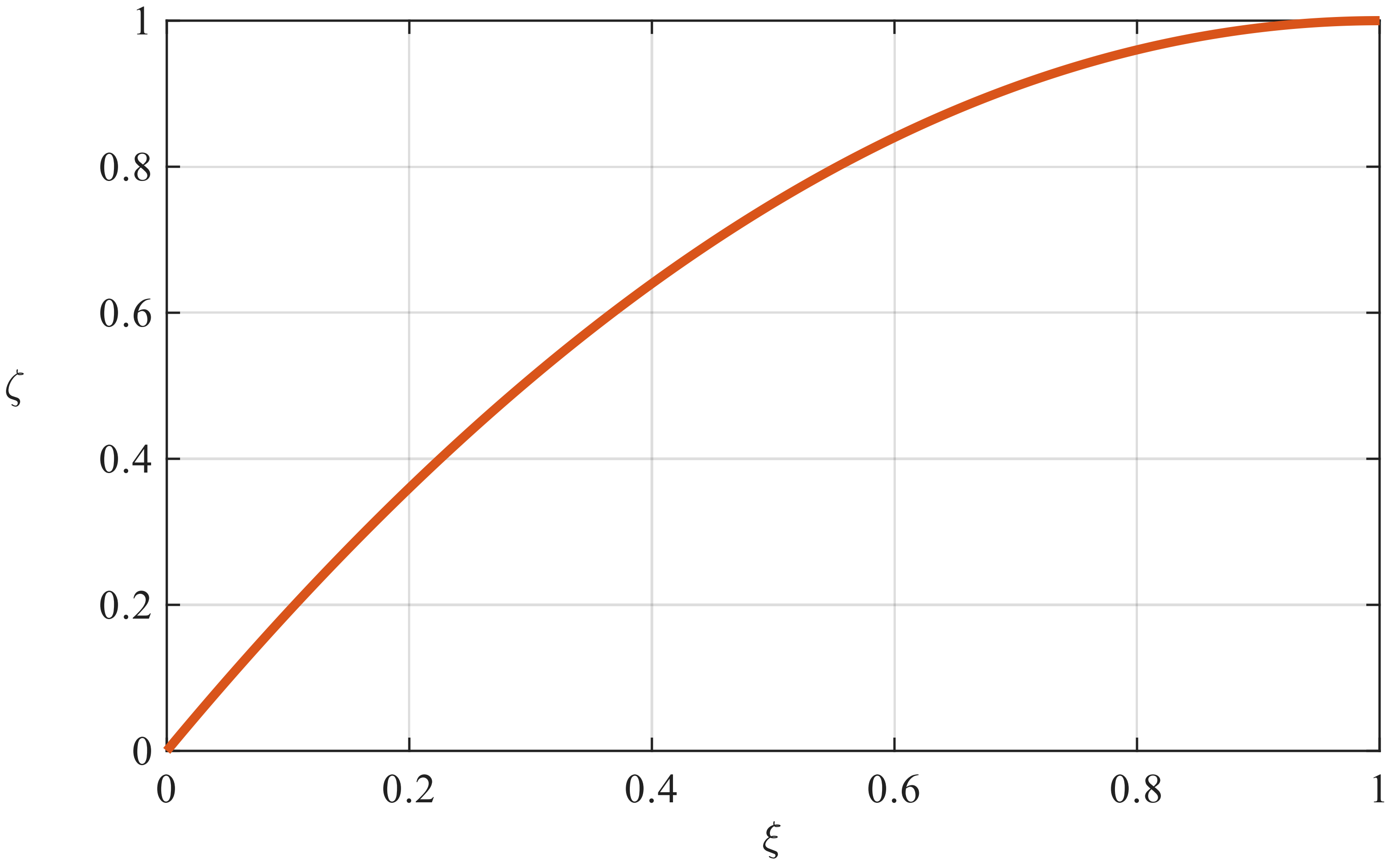}%
\caption{Relationship between perturbation amplitude and degree of steady-state breakdown}
\label{fig:3.1}
\end{figure}

For a gaseous star surrounded by vacuum undergoing gravitational collapse, the degree of steady-state breakdown characterizes how far the system departs from hydrostatic equilibrium, and it is an essential physical quantity for describing the collapse process. Using the perturbation amplitude parameter defined above, we define the degree of steady-state breakdown via a quadratic expression in the velocity ratio, namely
\begin{equation}\label{3.3}
\zeta (r)=1-\left(\frac {\mathcal{A} _{1}(r)}{v_{0}(r)}\right)^{2}=2\xi(r)-\xi^{2}(r).
\end{equation}
Moreover, there exists a direct physical connection between the perturbation amplitude parameter $\xi(r)$ and the degree of steady-state breakdown $\zeta(r)$ (see Figure~\ref{fig:3.1}), and the relationship between these two quantities is determined by the evolutionary properties of the system.

\subsection{Quantitative relationship between critical time $s_{1}(r)$ and perturbation amplitude parameter $\xi(r)$}
\indent\indent To better describe the quantitative relationships of the critical time $s_{1}(r)$ and collapse time $s_{0}(r)$ with the perturbation amplitude parameter $\xi(r)$ respectively, we introduce the following lemma.

\begin{lemma}\label{lem:3.1}
The relation between the initial velocity $\mathcal{A}_{1}(r)$ and $\xi(r)$ is given by
\begin{equation}\label{3.4}
\mathcal{A}_{1}(r)=\sqrt{2 G(r)}(1-\xi(r)).
\end{equation}
Moreover, the initial energy $n_{0}$ from \eqref{2.8} satisfies
\begin{equation}\label{3.5}
n_{0}=-2G(r)\xi(r) (2-\xi(r) ).
\end{equation}
\end{lemma}

\begin{proof}First, by \eqref{2.5} and \eqref{3.1}, we get
$$
\begin{aligned}
\mathcal{A} _{1}(r)&=v_{0}(r)(1-\xi(r))\\
&=\sqrt{2 G(r)}(1-\xi(r)),
\end{aligned}
$$
which completes the derivation of \eqref{3.4}.

Then again, by \eqref{2.5}, \eqref{2.8} and \eqref{3.1}, it yields
\[
\begin{aligned}
n_{0} & =\mathcal{A}_{1}^{2}-2 G(r)=2 G(r)\left(\frac{\mathcal{A}_{1}^{2}}{v_{0}^{2}(r)}-1\right) \\
& =2 G(r)\left((1-\xi(r))^{2}-1\right) \\
& =-2 G(r) \xi(r)(2-\xi(r)) .
\end{aligned}
\]
This verifies \eqref{3.5}, which completes the proof of \lemref{lem:3.1}.
\end{proof}

\begin{theorem}\label{thm:3.1}
For spherically symmetric collapse in the pressureless Euler-Poisson system, if $n_{0}<0$ and $0<\mathcal{A}_{1}(r)<v_{0}(r)$, the perturbation-based quantitative formula for the critical time $s_{1}(r)$ reads
\begin{align}\label{3.6}
s_{1}(r)=&\frac{1}{\sqrt{2 G(r)}(1-\xi(r))}\left(\frac{1}{\xi(r)(2-\xi(r))}-1\right) \nonumber \\
&+\frac{1}{\sqrt{2 G(r)}[\xi(r)(2-\xi(r))]^{\frac{3}{2}}} \arctan \left(\frac{1-\xi(r)}{\sqrt{\xi(r)(2-\xi(r))}}\right) .
\end{align}
\end{theorem}

\begin{proof}For $s<s_{1}$, we have $\mathcal{A}_{s}=\sqrt{\frac{2 G(r)}{\mathcal{A}}+n_{0}}$, which yields 
\begin{equation}\label{3.7}
\frac{d \mathcal{A}}{\sqrt{\frac{2 G(r)}{\mathcal{A}}+n_{0}}}=d s.
\end{equation}
Let $\eta=\sqrt{\frac{2 G(r)}{\mathcal{A}}+n_{0}}$, it follows that
\[
\frac{d \mathcal{A}}{\eta}=ds, \mathcal{A}=\frac{2G(r)}{{\eta}^{2}-n_{0}}, d \mathcal{A}=\frac{-4G(r)\eta}{(\eta^{2}-n_{0})^{2}}d\eta
\]
Integrating \eqref{3.7} with respect to $s$, we get
\begin{equation}\label{3.8}
\begin{aligned}
s=\int \frac{d \mathcal{A}}{\sqrt{\frac{2 G(r)}{\mathcal{A}}+n_{0}}}=&\int \frac{1}{\eta} \cdot\left(-\frac{4 G(r) \eta}{\left(\eta^{2}-n_{0}\right)^{2}}\right) d \eta=-4 G(r) \int \frac{1}{\left(\eta^{2}-n_{0}\right)^{2}} d \eta\\
=&\frac{2 G(r)}{\eta \left(\eta^{2}-n_{0}\right)}+\frac{2 G(r)}{n_{0}} \cdot \frac{1}{\eta}+\frac{2 G(r)}{n_{0} \sqrt{-n_{0}}} \arctan \left(\frac{\eta}{\sqrt{-n_{0}}}\right)+C_{1} \\
=&\frac{2 G(r) \eta}{n_{0}\left(\eta^{2}-n_{0}\right)}+\frac{2 G(r)}{n_{0} \sqrt{-n_{0}}} \arctan \left(\frac{\eta}{\sqrt{-n_{0}}}\right)+C_{1} .
\end{aligned}
\end{equation}
Since $\mathcal{A}(0,r)=1$, it yields
\begin{equation}\label{3.9}
\eta(0)=\sqrt{2 G(r)+n_{0}}=\sqrt{2 G(r)+\mathcal{A}_{1}^{2}-2 G(r)}=\mathcal{A}_{1}.
\end{equation}
Substituting \eqref{3.9} into \eqref{3.8}, we get
\begin{equation}\label{3.10}
C_{1}=-\left[\frac{1}{\mathcal{A}_{1}}+\frac{2 G(r)}{n_{0} \mathcal{A}_{1}} +\frac{2 G(r)}{n_{0} \sqrt{-n_{0}}} \arctan \left(\frac{\mathcal{A}_{1}}{\sqrt{-n_{0}}}\right)\right].
\end{equation}
At the critical time $s=s_{1}$, we have $\eta=0$. Substitute $\eta=0$ into \eqref{3.8}, we get
\begin{equation}\label{3.11}
s_{1}(r)=-\left[\frac{1}{\mathcal{A}_{1}}+\frac{2 G(r)}{n_{0} \mathcal{A}_{1}} + \frac{2 G(r)}{n_{0} \sqrt{-n_{0}}} \arctan \left(\frac{\mathcal{A}_{1}}{\sqrt{-n_{0}}}\right)\right].
\end{equation}
From \lemref{lem:3.1}, we know
\[
\sqrt{-n_{0}}=\sqrt{2 G(r) \xi(r)(2-\xi(r))},
\]
which gives the ratio
\[
\frac{\mathcal{A}_{1}}{\sqrt{-n_{0}}}=\frac{1-\xi(r)}{\sqrt{\xi(r)(2-\xi(r))}}, \quad \frac{1}{\mathcal{A}_{1}(r)}=\frac{1}{\sqrt{2 G(r)}(1-\xi(r))},
\]
and 
\[
\frac{2 G(r)}{n_{0} \mathcal{A}_{1}(r)}=-\frac{1}{\sqrt{2 G(r)} \xi(r)(1-\xi(r))(2-\xi(r))} .
\]
Therefore, \eqref{3.11} becomes
\[
\begin{aligned}
s_{1}(r)=&\frac{1}{\sqrt{2 G(r)}(1-\xi(r))}\left(\frac{1}{\xi(r)(2-\xi(r))}-1\right)\\
+&\frac{1}{\sqrt{2 G(r)}[\xi(r)(2-\xi(r))]^{\frac{3}{2}}} \arctan \left(\frac{1-\xi(r)}{\sqrt{\xi(r)(2-\xi(r))}}\right) .
\end{aligned}
\]
This completes the proof of \theref{thm:3.1}.
\end{proof}

\begin{remark}
From \eqref{3.3} and \theref{thm:3.1}, the relationship between the degree of steady-state breakdown and the critical time follows directly, and we omit detailed elaboration of this relation in the present paper.
\end{remark}

This conclusion reveals the monotonic nature of gravitational collapse for gaseous stars. It can be utilized to infer relevant physical properties throughout the collapse process, provide support for subsequent quantitative research on collapse dynamics, and deepen our fundamental understanding of the essential mechanism governing gravitational collapse in gaseous stars.

\subsection{Quantitative relationship between collapse time $s_{0}(r)$ and perturbation amplitude parameter $\xi(r)$}
\indent\indent This subsection derives and proves the analytical formula of collapse time \(s_0(r)\) with perturbation amplitude \(\xi(r)\). Asymptotic analysis yields the power-law relation \(s_0(r)\propto\xi(r)^{-\frac{3}{2}}\) for weak perturbations. Theoretical results show larger \(\xi(r)\) speeds up collapse, and collapse time grows with radial coordinate at fixed perturbation, quantitatively describing the gravitational collapse of gaseous stars.

\begin{theorem}\label{thm:3.2}
For spherically symmetric collapse within the pressureless Euler-Poisson system, under the conditions $n_{0}<0$ and $0<\mathcal{A}_{1}(r)<v_{0}(r)$, the perturbation-based quantitative formula for the collapse time $s_{0}(r)$ is given by
\begin{align}\label{3.12}
s_{0}(r)&=\frac{1}{\sqrt{2G(r)}(1-\xi(r))}\left(\frac{1}{\xi(r)(2-\xi(r))}-1\right) \nonumber \\
&+\frac{1}{\sqrt{2 G(r)}[\xi(r)(2-\xi(r))]^{\frac{3}{2}}} \arctan \left(\frac{1-\xi(r)}{\sqrt{\xi(r)(2-\xi(r))}}\right) \nonumber \\
&+\frac{\pi}{2 \sqrt{2 G(r)} \cdot[\xi(r)(2-\xi(r))]^{\frac{3}{2}}}.
\end{align}
\end{theorem}

\begin{proof}For $s>s_{1}$, we have $\mathcal{A}_{s}=-\sqrt{\frac{2 G(r)}{\mathcal{A}}+n_{0}}=-\eta<0$, which gives
\begin{equation}\label{3.13}
\frac{d \mathcal{A}}{-\sqrt{\frac{2 G(r)}{\mathcal{A}}+n_{0}}}=d s.
\end{equation}
Integrating both sides yields
\begin{equation}\label{3.14}
s=-\frac{2 G(r)}{\eta \left(\eta^{2}-n_{0}\right)}-\frac{2 G(r)}{n_{0}} \cdot \frac{1}{\eta}-\frac{2 G(r)}{n_{0} \sqrt{-n_{0}}} \arctan \left(\frac{\eta}{\sqrt{-n_{0}}}\right)+s_{1}.
\end{equation}
As $s \to s_{0}$, $\mathcal{A} \to 0$ and $\eta \to +\infty$. Then it follows from \eqref{3.14}
\begin{equation}\label{3.15}
s_{0}=s_{1}-\frac{G(r) \pi}{n_{0} \sqrt{-n_{0}}}.
\end{equation}
From \theref{thm:3.1} and the identity
\[
\frac{G(r) \pi}{n_{0} \sqrt{-n_{0}}}=-\frac{\pi}{2 \sqrt{2 G(r)} \cdot[\xi(r)(2-\xi(r))]^{\frac{3}{2}}},
\]
we deduce
\begin{align}
s_{0}(r)&=\frac{1}{\sqrt{2G(r)}(1-\xi(r))}\left(\frac{1}{\xi(r)(2-\xi(r))}-1\right) \nonumber \\
&+\frac{1}{\sqrt{2 G(r)}[\xi(r)(2-\xi(r))]^{\frac{3}{2}}} \arctan \left(\frac{1-\xi(r)}{\sqrt{\xi(r)(2-\xi(r))}}\right) \nonumber \\
&+\frac{\pi}{2 \sqrt{2 G(r)} \cdot[\xi(r)(2-\xi(r))]^{\frac{3}{2}}} \nonumber.
\end{align}
This completes the proof of \theref{thm:3.2}.
\end{proof}

\begin{remark}
From \eqref{3.3} and \theref{thm:3.2}, the relationship between the degree of steady-state breakdown and the collapse time follows directly, and we omit further detailed elaboration in this paper.
\end{remark}

\begin{proposition}\label{prop:3.1}
Let $0<\xi(r)<1$, so that $n_{0}<0$, and set
$\lambda=\xi(r)(2-\xi(r))$. Then, as $\xi(r)\to0^{+}$,
\begin{equation}\label{3.16}
s_{0}(r)=\frac{\pi}{4\sqrt{G(r)}}\,\xi(r)^{-\frac{3}{2}}
\bigl(1+O(\xi(r))\bigr),
\end{equation}
and consequently $s_{0}(r)\propto \xi(r)^{-3/2}$. The estimate is uniform for $r$
in any compact subinterval of $(0,1]$ on which $G$ is bounded away from zero.
\end{proposition}

\begin{proof}
For any fix $r\in(0,1]$, we write $G=G(r)$, $\xi=\xi(r)$, $\lambda=\xi(2-\xi)$, then $n_{0}=-2G\lambda$ and $1-\xi=\sqrt{1-\lambda}$. By Theorem~\ref{thm:3.1} and \eqref{3.15}, it yields
\begin{equation}\label{3.17}
s_{0}=s_{1}-\frac{G\pi}{n_{0}\sqrt{-n_{0}}}
=\frac{1}{\sqrt{2G}}\left[
\frac{1}{1-\xi}\Bigl(\frac{1}{\lambda}-1\Bigr)
+\frac{1}{\lambda^{3/2}}\Bigl(\arctan\frac{1-\xi}{\sqrt{\lambda}}+\frac{\pi}{2}\Bigr)
\right].
\end{equation}
Since $\arctan\sqrt{\lambda^{-1}-1}+\frac{\pi}{2}=\pi-\arcsin\sqrt{\lambda}$ and $1-\xi=\sqrt{1-\lambda}$, \eqref{3.17} becomes
\begin{equation}\label{3.18}
s_{0}=\frac{1}{\sqrt{2G}}\left[
\frac{\sqrt{1-\lambda}}{\lambda}+\frac{\pi-\arcsin\sqrt{\lambda}}{\lambda^{3/2}}
\right].
\end{equation}
Using $\sqrt{1-\lambda}=1-\frac{\lambda}{2}+O(\lambda^{2})$ and $\arcsin\sqrt{\lambda}=\sqrt{\lambda}+\frac{1}{6}\lambda^{3/2}+O(\lambda^{5/2})$, we get
\[
\frac{\sqrt{1-\lambda}}{\lambda}=\frac{1}{\lambda}-\frac{1}{2}+O(\lambda),
\quad
\frac{\pi-\arcsin\sqrt{\lambda}}{\lambda^{3/2}}
=\frac{\pi}{\lambda^{3/2}}-\frac{1}{\lambda}-\frac{1}{6}+O(\lambda).
\]
So
\begin{equation}\label{3.19}
s_{0}=\frac{1}{\sqrt{2G}}\left[
\frac{\pi}{\lambda^{3/2}}-\frac{2}{3}+O(\lambda)
\right],\quad \lambda\to0^{+}.
\end{equation}
$\lambda=2\xi\left(1-\frac{\xi}{2}\right)$ gives $\lambda^{-3/2}=2^{-3/2}\xi^{-3/2}\left(1+\frac{3}{4}\xi+O(\xi^{2})\right)$. Substituting it into \eqref{3.19}, we obtain the refined expansion
\begin{equation}\label{3.20}
s_{0}=\frac{\pi}{4\sqrt{G}}\,\xi^{-3/2}
+\frac{3\pi}{16\sqrt{G}}\,\xi^{-1/2}
-\frac{2}{3\sqrt{2G}}+O(\xi^{1/2}).
\end{equation}
Since all non-leading terms are $O(\xi^{-1/2})$, this gives 
\begin{equation}\label{3.21}
s_{0}=\frac{\pi}{4\sqrt{G}}\,\xi^{-3/2}\bigl(1+O(\xi)\bigr),
\end{equation}
Here the $O(\xi)$ in \eqref{3.21} is relative, i.e. it multiplies the leading term, \eqref{3.21} is also equivalent to $s_0=\frac{\pi}{4\sqrt G}\xi^{-3/2}+O(\xi^{-1/2})$. Thus we proved \eqref{3.16}.
\end{proof}

\begin{corollary}\label{cor:3.1}
Under the assumptions of Proposition~\ref{prop:3.1}, the following hold.
\begin{enumerate}
\item[(i)] The relative error of the asymptotic law satisfies
\[
R(\xi):=\left|\frac{s_{0}(r)}{\frac{\pi}{4\sqrt{G(r)}}\,\xi^{-\frac{3}{2}}}-1\right|
=\frac{3}{4}\,\xi+O(\xi^{\frac{3}{2}}),\quad \xi\to0^{+}.
\]
\item[(ii)] The two time scales satisfy
\[
\lim_{\xi\to0^{+}}\frac{s_{0}}{s_{1}}=2,
\quad
\lim_{\xi\to1^{-}}\frac{s_{0}}{s_{1}}=+\infty .
\]
\end{enumerate}
\end{corollary}

\begin{proof}
Part (i) follows \eqref{3.20}. Setting $A=\frac{\pi}{4\sqrt{G}}\xi^{-3/2}$, we have
\[
s_{0}=A\left(1+\frac{3}{4}\xi+O(\xi^{3/2})\right),
\]
and hence $R(\xi)=\frac34\xi+O(\xi^{3/2})$.

For (ii), write $\lambda=\xi(2-\xi)$ and $\mathcal A_1=\sqrt{2G}(1-\xi)$.
By \eqref{3.6}, with $a=\sqrt{-n_0}=\sqrt{2G\lambda}$,
\[
s_1=\frac{1}{\lambda\mathcal A_1}
+\frac{1}{\lambda\sqrt{2G\lambda}}\arctan\frac{\mathcal A_1}{a}+O(1).
\]
Since $\lambda=2\xi(1-\xi/2)$, we have $\sqrt{\lambda}=\sqrt{2\xi}(1+O(\xi))$;
hence $\mathcal A_1/a=(1-\xi)/\sqrt{\lambda}=(1+O(\xi))/\sqrt{2\xi}$, and then
$\arctan(\mathcal A_1/a)=\frac{\pi}{2}-O(\xi^{1/2})$. Thus
\[
s_1=\frac{\pi}{2\sqrt{2G}}\,\lambda^{-3/2}\bigl(1+O(\xi^{1/2})\bigr).
\]
Set
\[
P=\frac{\pi}{2\sqrt{2G}}\,\lambda^{-3/2}.
\]
The preceding estimate gives
\[
s_1=P\bigl(1+O(\xi^{1/2})\bigr)=P+O(\xi^{-1}).
\]
By \eqref{3.15}, we get \(s_0=s_1+P\), and hence
\[
s_0=2P+O(\xi^{-1}),
\quad
2s_1=2P+O(\xi^{-1}).
\]
Consequently, \(s_0-2s_1=O(\xi^{-1})\). Since $s_1\sim P\sim\xi^{-3/2}$, it yields
\[
\frac{s_0-2s_1}{s_1}=O(\xi^{1/2}),
\]
then 
\[
\quad
\frac{s_0}{s_1}=2\bigl(1+O(\xi^{1/2})\bigr)\longrightarrow2, \quad as \quad \xi \to 0^{+}.
\]
As $\xi\to1^{-}$, from formula \eqref{3.6} and \eqref{3.12}, we get $s_{1}\to0$ and $s_{0}\to\frac{\pi}{2\sqrt{2G}}$. Hence $\lim\limits_{\xi\to1^{-}}\frac{s_0}{s_1}=+\infty$.
\end{proof}

The study on the relationship between $\xi(r)$ and $s_{0}(r)$ reveals that the system is highly sensitive to variations in perturbation amplitude. As $\xi(r)$ increases gradually, the collapse time decreases rapidly and eventually approaches zero. This provides a quantitative description of the collapse behaviors of gaseous stars under perturbations of different magnitudes.

\begin{figure}[H]
\centering
  \includegraphics[width=0.8\textwidth]{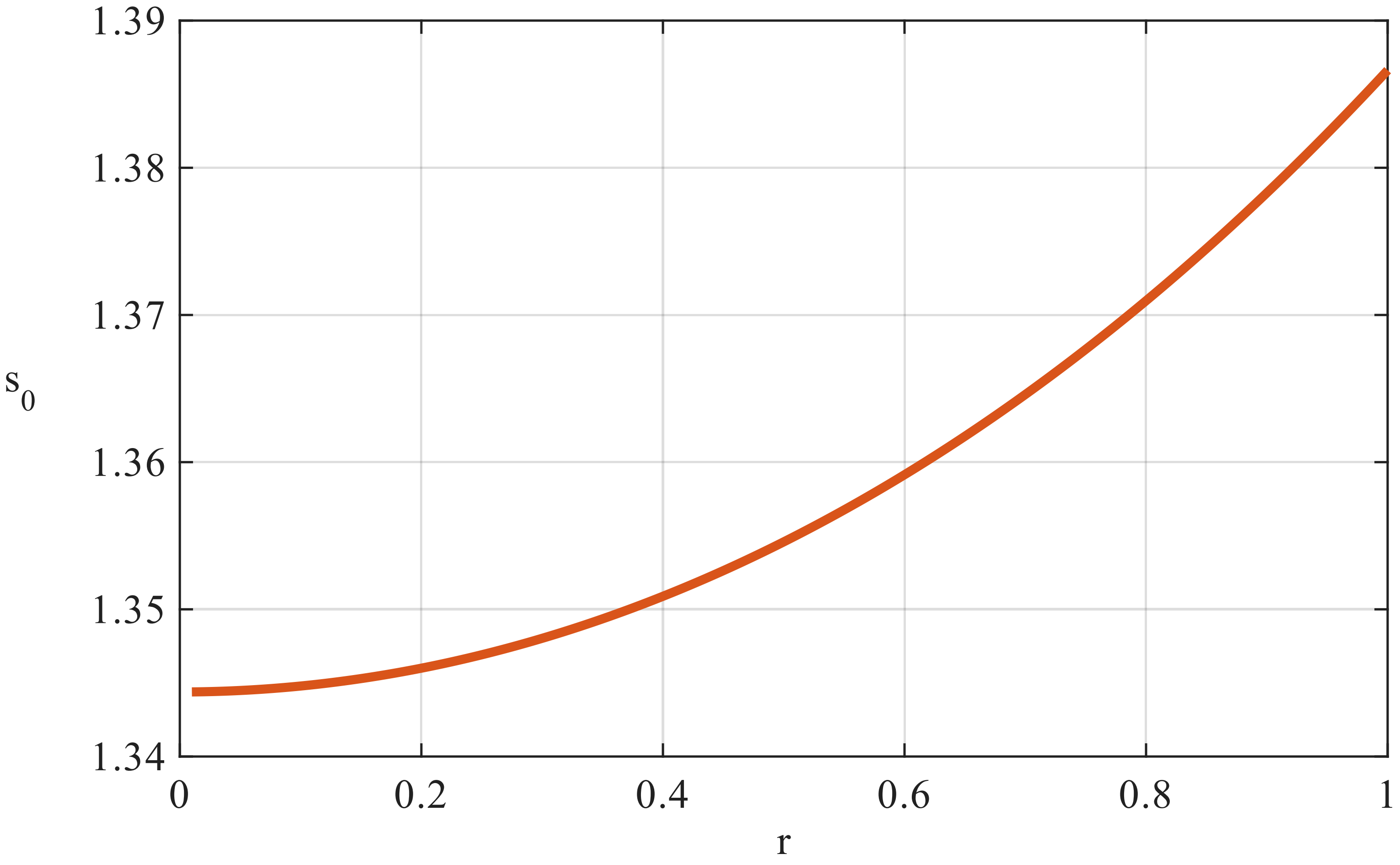}%
\caption{Variation of collapse time with radial coordinate for fixed perturbation amplitude}
\label{fig:3.2}
\end{figure}

Meanwhile, it is easy to observe that when the perturbation amplitude is fixed, the collapse time is positively correlated with the radial coordinate (see Figure~\ref{fig:3.2}).

\section{Numerical verification and quantitative behavior}

\indent\indent This section has three purposes. First, we verify the closed-form formulas
\eqref{3.6} and \eqref{3.12} against an independent numerical integration of \eqref{3.7}. Second, we
quantify the accuracy of the asymptotic law of Proposition~\ref{prop:3.1}. Third, we
compare the two time scales $s_1$ and $s_0$ and examine their ratio. Since the existence
and the qualitative classification of collapse solutions were established in
\cite{huang2022}, we do not reproduce the trajectories here; the emphasis is on the
quantitative relations among $\xi$, $\zeta$, $s_1$ and $s_0$.

\subsection{Verification and behavior of the closed-form formulas}

\indent\indent All computations are carried out at the representative radius $r=0.5$,
for which $G(r)=4.126$ and $v_0(r)=2.873$. For a given $\xi$, the numerical values
$s_1^{\mathrm{num}}$ and $s_0^{\mathrm{num}}$ are obtained by integrating \eqref{3.7}
directly. The integrand has an integrable singularity at $\mathcal A=\mathcal A_{\max}$;
the substitution $\mathcal A=\mathcal A_{\max}-u^2$ removes it, and the resulting smooth
integral is evaluated by adaptive quadrature with tolerance $10^{-12}$. The numerical
values at $\xi=0.01,0.05,0.10,0.30,0.60,0.90$ reproduce the closed-form curves \eqref{3.6}
and \eqref{3.12} to within a maximum relative deviation of $3\times10^{-9}$.At each shell $\xi=\xi(r)$. Both the closed forms \eqref{3.6} and \eqref{3.12} and the
integral in \eqref{3.7} factor as $G(r)^{-\frac{1}{2}}F(\xi)$, with the same function $F$ for
every shell; the only other occurrence of $r$ is through the value $\xi(r)$ itself, the upper limit $\mathcal A_{\max}=\frac{1}{\xi(2-\xi)}$ being independent of $G$. Hence
the relative deviation is a function of $\xi$ alone, and the scan over $\xi$ at
$r=0.5$ represents all shells.

\begin{figure}[H]
\centering
\includegraphics[width=0.70\textwidth]{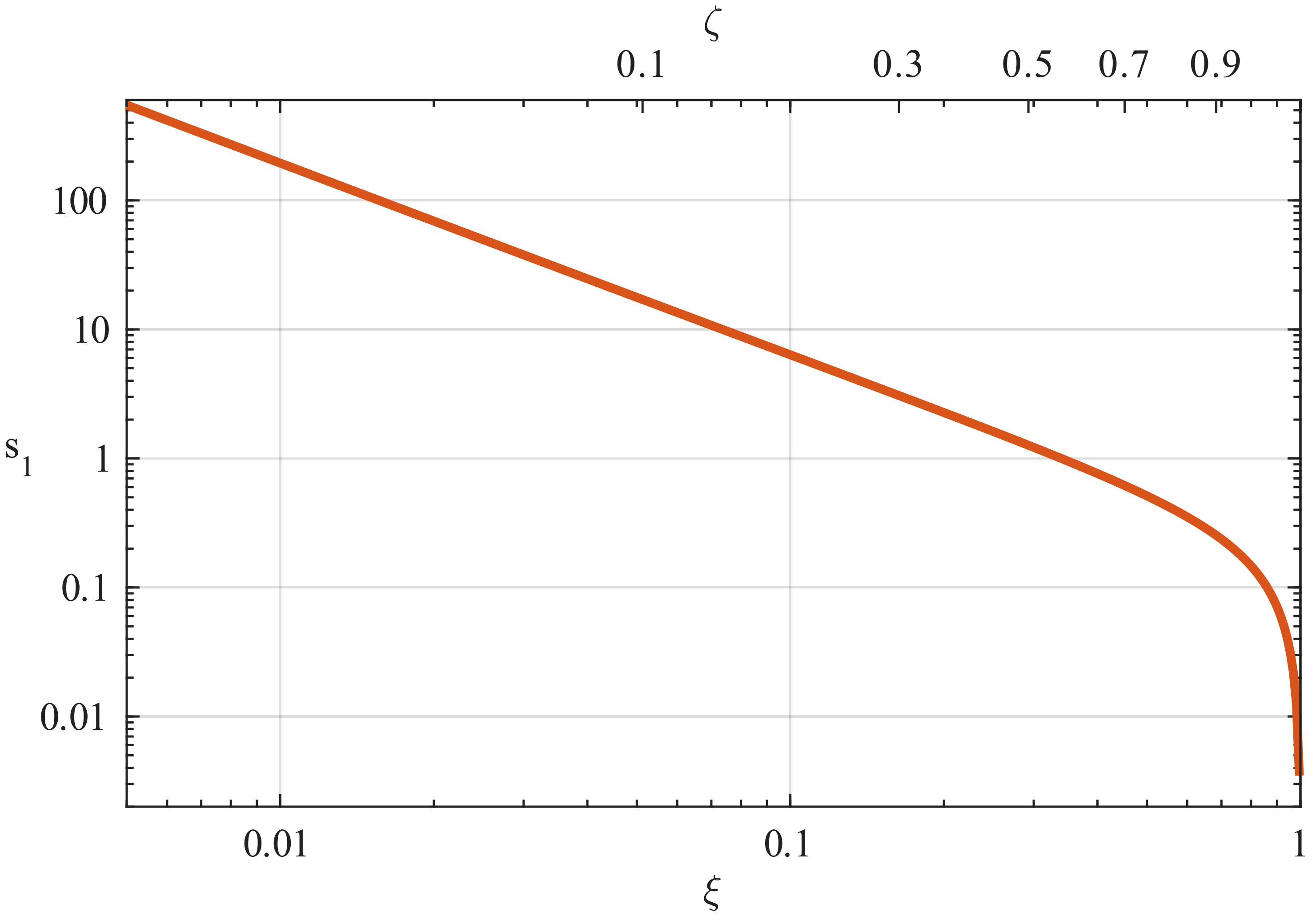}%
\caption{Critical time $s_1$ as a function of the perturbation amplitude $\xi$; the
upper axis gives the steady-state breakdown degree $\zeta$.}
\label{fig:4.1}
\end{figure}

\begin{figure}[H]
\centering
\includegraphics[width=0.70\textwidth]{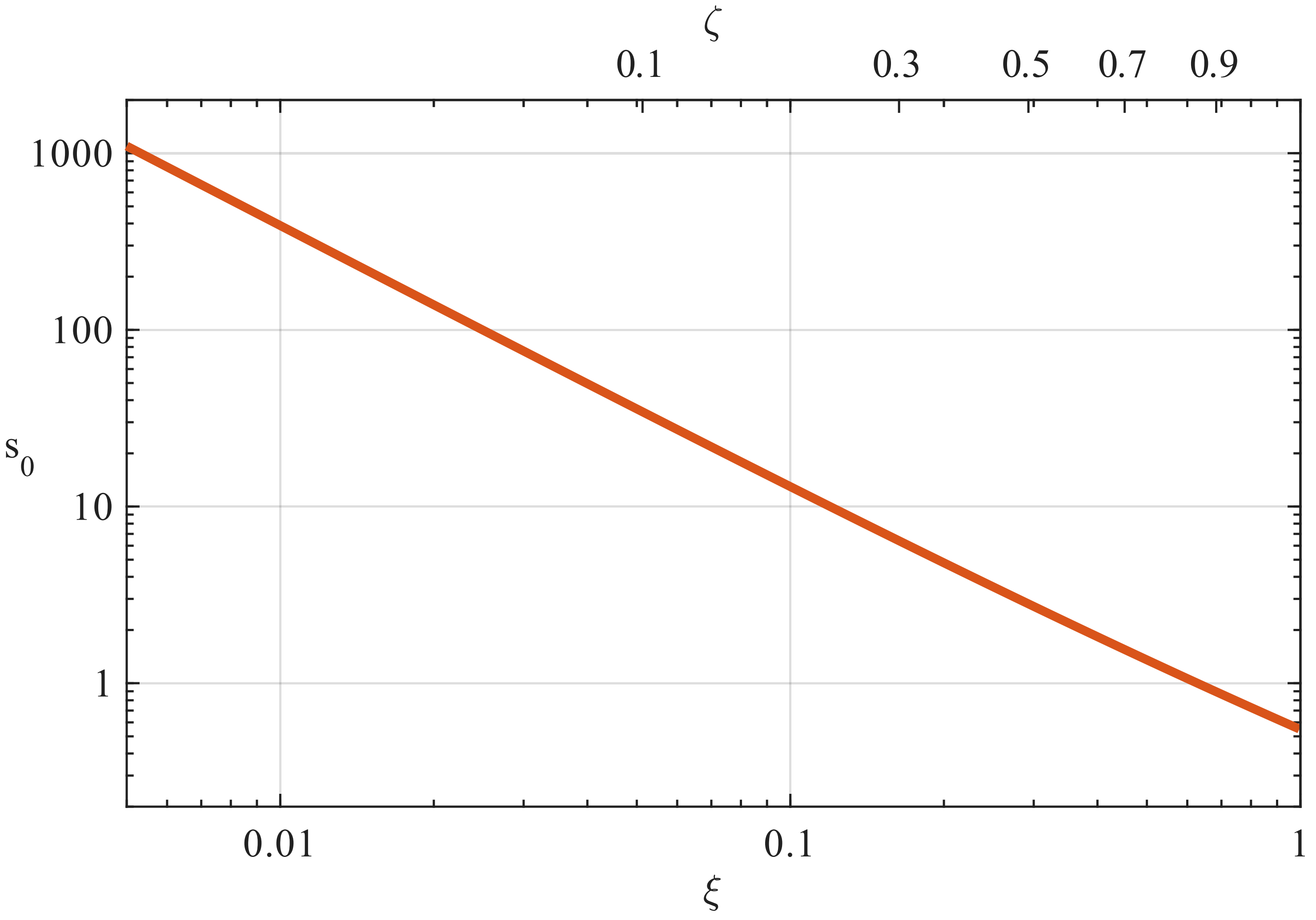}%
\caption{Collapse time $s_0$ as a function of the perturbation amplitude $\xi$; the
upper axis gives the steady-state breakdown degree $\zeta$.}
\label{fig:4.2}
\end{figure}

Figures~\ref{fig:4.1} and~\ref{fig:4.2} display the two closed-form curves together
with the $\zeta$ axis. Both $s_1$ and $s_0$ decrease strictly with $\xi$; they diverge
as $\xi\to0^{+}$, whereas $s_1\to0$ and $s_0\to\frac{\pi}{2\sqrt{2G}}$ as $\xi\to1^{-}$, the
latter being the free-fall time from rest. The two curves have the same qualitative
shape and differ essentially by the factor $\frac{s_0}{s_1}\in[2,+\infty)$ discussed below.

\subsection{Asymptotic validity and the two time scales}

\begin{figure}[H]
\centering
  \includegraphics[width=0.70\textwidth]{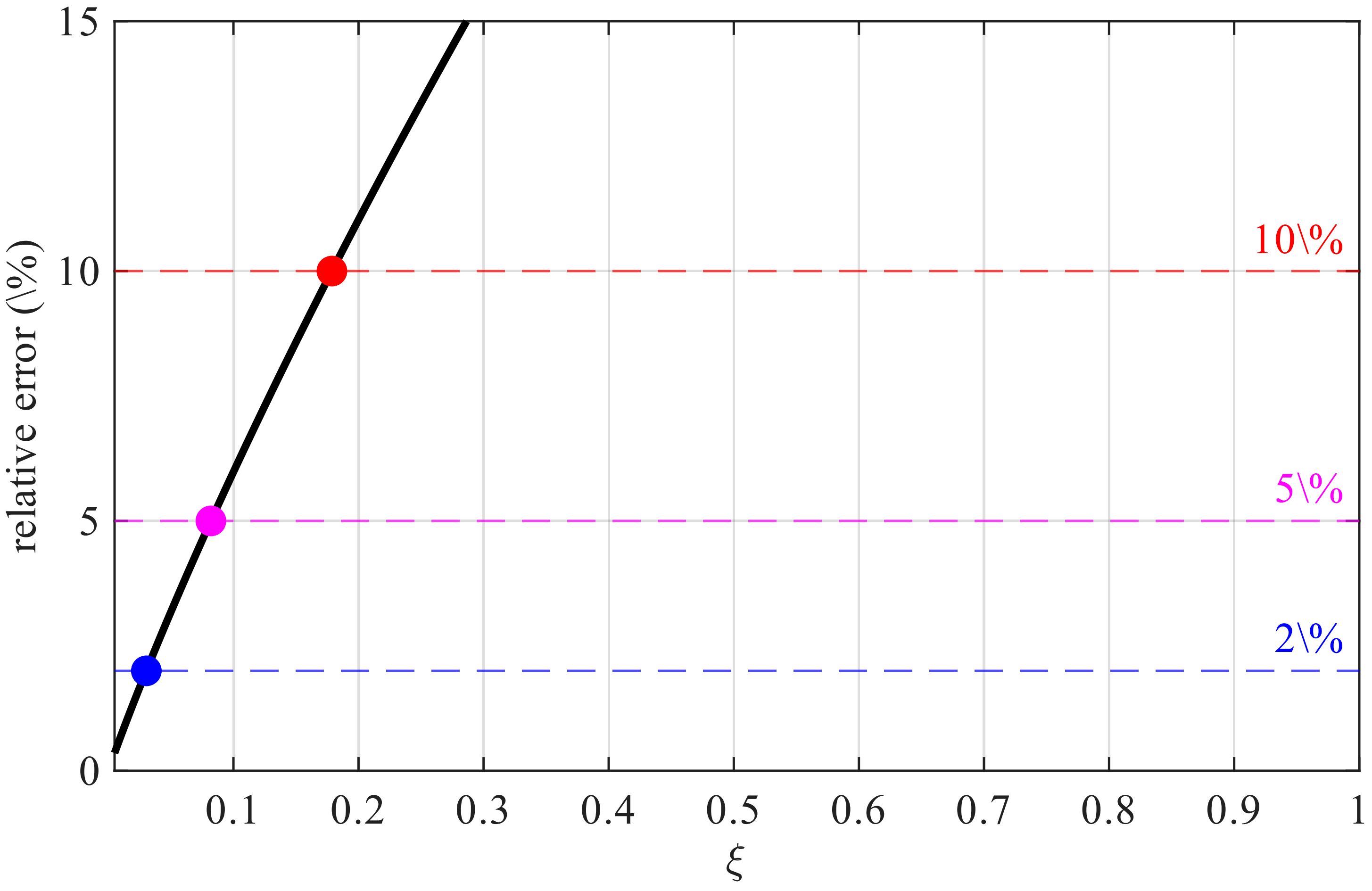}%
\caption{Relative error of the asymptotic law as a function of $\xi$; the dashed lines
mark the $2\%$, $5\%$ and $10\%$ levels.}
\label{fig:4.3}
\end{figure}

\indent\indent Figure~\ref{fig:4.3} shows the relative error $R(\xi)$ of the asymptotic law defined in
Corollary~\ref{cor:3.1}(i), which indicates $R(\xi)=\frac34\xi+O(\xi^{\frac{3}{2}})$ as
$\xi\to0^{+}$. The computed values give $R<2\%$ for $\xi\le0.03$, $R<5\%$ for
$\xi\le0.09$ and $R<10\%$ for $\xi\le0.20$, whereas $R(0.30)=13.5\%$. In terms of
$\zeta=\xi(2-\xi)$, these thresholds correspond to $\zeta\approx0.06$, $0.17$ and
$0.36$. Accordingly, the power law is accurate to $2\%$ for $\xi\le0.03$ and to $10\%$
for $\xi\le0.2$; beyond $\xi\approx0.3$, where the relative error exceeds $13\%$, it
should be used only as an order-of-magnitude estimate.

\begin{figure}[H]
\centering
  \includegraphics[width=0.70\textwidth]{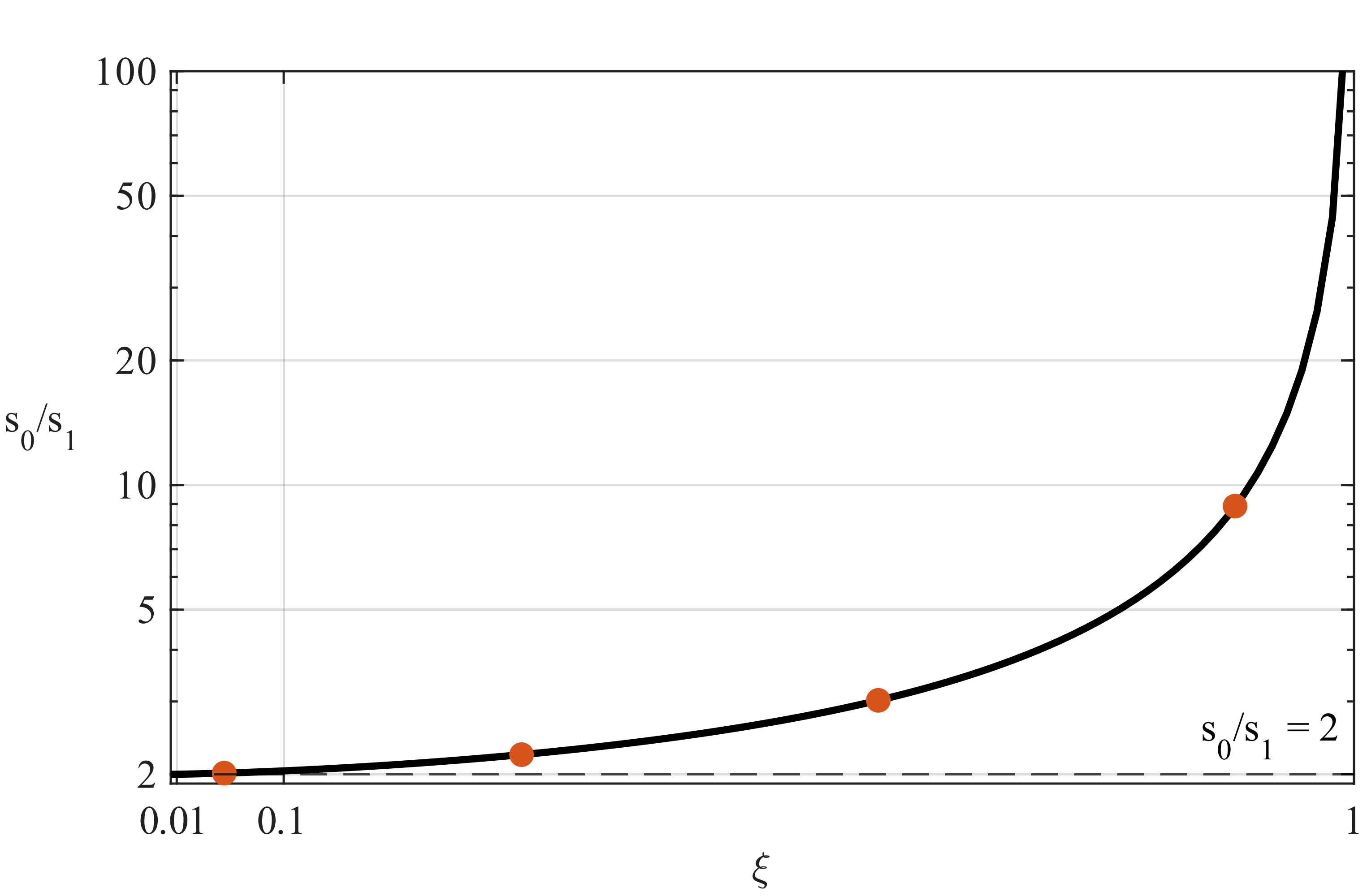}%
\caption{Ratio $s_0/s_1$ as a function of $\xi$; the dashed line marks the limit $2$
as $\xi\to0^{+}$.}
\label{fig:4.4}
\end{figure}

The two time scales are compared in Figure~\ref{fig:4.4}. By
Corollary~\ref{cor:3.1}(ii), we get $\frac{s_0}{s_1}\to2$ as $\xi\to0^{+}$: for weak perturbations,
the expansion and collapse phases last about equally long, and already at $\xi=0.05$
the ratio is $2.01$. As $\xi\to1^{-}$ the expansion phase disappears, $s_1\to0$ while
$s_0\to\frac{\pi}{2\sqrt{2G}}$, so the ratio diverges; numerically $\frac{s_0}{s_1}=2.23$ at
$\xi=0.3$, $3.02$ at $\xi=0.6$ and $8.87$ at $\xi=0.9$.

\section{Conclusion}

\indent\indent This paper has studied the gravitational collapse of gaseous stars governed
by the pressureless Euler-Poisson system in the regime $1<\gamma<4/3$. Main
contribution is a set of explicit quantitative relations linking the initial velocity
deficit to the critical and collapse times.

Introducing the perturbation amplitude $\xi(r)=1-\frac{\mathcal A_1(r)}{v_0(r)}$ and the
breakdown measure $\zeta(r)=2\xi-\xi^2$, we derived the closed-form expressions \eqref{3.6}
and \eqref{3.12} for $s_1(r)$ and $s_0(r)$ at each Lagrangian radius. These yield the
asymptotic law
\[
s_0(r)=\frac{\pi}{4\sqrt{G(r)}}\,\xi(r)^{-\frac{3}{2}}\bigl(1+O(\xi(r))\bigr), \quad as \quad \xi(r)\to0^{+},
\]
together with the relative-error expansion $R(\xi)=\frac34\xi+O(\xi^{\frac{3}{2}})$; the
exponent $-\frac{3}{2}$ is independent of $r$, while the amplitude is fixed by $G(r)$. The same
relations give the two-time-scale limit $s_0/s_1\to2$ as $\xi\to0^{+}$ and divergence as
$\xi\to1^{-}$.

Numerical integration of \eqref{3.7} reproduces \eqref{3.6} and \eqref{3.12} to better than $10^{-8}$.
The asymptotic law is accurate to $2\%$ for $\xi\lesssim0.03$, $5.6\%$ at $\xi=0.1$ and
$13.5\%$ at $\xi=0.3$; with a $10\%$ tolerance the accuracy threshold is
$\xi\approx0.2$ (equivalently $\zeta\approx0.36$). These results provide a simple
quantitative estimate of the collapse timescale and clarify the two-phase structure of
the collapse.

Two directions remain open: extending the quantitative analysis to the regimes
$\gamma\ge4/3$, which lie outside the pressureless-collapse setting considered here,
and describing stronger perturbations $\xi\gtrsim0.3$, for which higher-order terms of
the expansion or the full Euler-Poisson system must be retained.

\end{sloppypar}
\end{document}